\documentclass{article}
\usepackage[T2A]{fontenc}
\usepackage[cp1251]{inputenc}
\usepackage[english,russian]{babel}
\usepackage[tbtags]{amsmath}
\usepackage{amsfonts,amssymb,mathrsfs,amscd}
\usepackage{verbatim}
\usepackage{eucal}
\usepackage{graphicx}
\usepackage{enumerate}
\usepackage[pdftex,unicode,colorlinks,linkcolor=blue,
citecolor=red,bookmarksopen,pdfhighlight=/N]{hyperref}
\usepackage[hyper]{faa-1-arxiv}

\JournalName{%%Функциональный анализ и его приложения
}
\numberwithin{equation}{section}
\theoremstyle{plain}
\newtheorem{theorem}{Теорема}
\newtheorem{theoA}{\bf Теорема A}
\newtheorem{theoB}{\bf Теорема B}
\newtheorem{lemma}{\bf Лемма}
\newtheorem{propos}{\bf Предложение}
\newtheorem{corollary}{\bf Следствие}
\theoremstyle{definition}
\newtheorem{definition}{Определение}

\newtheorem{proof}{Доказательство}
\newtheorem{example}{Пример}
\newtheorem{remark}{Замечание}
\newcommand{\const}{{\rm const}}

\renewcommand{\leq}{\leqslant}
\renewcommand{\geq}{\geqslant}
\renewcommand{\d}{\,{\rm d}}
\newcommand{\dd}{\,{\rm d}}
\newcommand{\e}{\varepsilon}
\def\RR{\mathbb R}
\def\CC{\mathbb C}
\def\NN{\mathbb N}
\def\DD{\mathbb D}

\DeclareMathOperator{\clos}{clos}
\DeclareMathOperator{\Int}{int}
\DeclareMathOperator{\Har}{har}
\DeclareMathOperator{\Hol}{Hol}
\DeclareMathOperator{\dist}{dist}
\DeclareMathOperator{\Zero}{Zero}
\DeclareMathOperator{\sbh}{sbh}
\DeclareMathOperator{\dom}{dom}
\DeclareMathOperator{\dsbh}{\text{$\delta${\rm -sbh}}}
\DeclareMathOperator{\supp}{supp}
\DeclareMathOperator{\comp}{c}
\DeclareMathOperator{\Meas}{Meas}
\DeclareMathOperator{\loc}{loc}
\begin{document}

\title{К распределению нулевых множеств голоморфных функций. III. Теоремы обращения}

\author[B.\,N.~Khabibullin]{Б.\,Н.~Хабибуллин}
%%\address{Башкирский государственный университет}
%%\email{khabib-bulat@mail.ru}
%%второй автор
\author[F.\,B.~Khabibullin]{Ф.\,Б.~Хабибуллин}
%%\address{Башкирский государственный университет}
%%\email{khabibullinfb@mail.ru}

%%%\date{01.11.2018}
\udk{517.53 : 517.574 : 517.987.1}

\maketitle

\begin{abstract}
Пусть   $M$ --- субгармоническая функции в области $D\subset  \CC^n$ с мерой Рисса $\nu_M$, ${\sf Z}\subset D$. Как было показано в первой из предшествующих статей, если существует голоморфная функция $f\neq 0$ в $D$, $f({\sf Z})=0$,  $|f|\leq \exp M$ на $D$, то имеет место некоторая {\it шкала\/} интегральных равномерных оценок сверху распределения множества $\sf Z$ через $\nu_M$.
В настоящей статье показано, что при $n=1$ этот результат <<почти обратим>>.  Из такой {\it шкалы\/} оценок распределения точек последовательности ${\sf Z}:=\{{\sf z}_k \}_{k=1,2,\dots}\subset D\subset \CC$ через $\nu_M$  следует, что существует  ненулевая голоморфной функции $f$ в $D$, $f(\sf Z)=0$, $|f|\leq \exp  M^{\uparrow}$ на $D$, где функция $M^{\uparrow}\geq M$ на $D$ строится через усреднения функции $M $ по  быстро сужающимся кругам при приближении к границе области $D$ с некоторой возможной  аддитивной логарифмической добавкой, связанной со скоростью сужения этих кругов.
\end{abstract}
	
\markright{К распределению нулевых множеств голоморфных функций. III. \dots}

\footnotetext[0]{Исследование выполнено за счёт гранта Российского
	научного фонда (проект № 18-11-00002 --- первый автор), а также при поддержке РФФИ 
(проекты №№ 16-01-00024, 18-51-06002 --- второй автор)}

\section{Введение}\label{In1}

В настоящей статье используются обозначения, определения,  соглашения и результаты из \cite{KhaRoz18} с их естественными адаптациями   для {\it комплексной плоскости\/} $\CC$ и ее {\it одноточечной компактификации Александрова\/}  $\CC_{\infty}:=\CC\cup \{\infty\}$. Основная цель  --- дать обращения результатов исходной работы \cite{KhaRoz18}  для областей $D\subset  \CC_{\infty}$ в форме критерия в субгармонической версии, а также в виде, близком к критерию, когда рассматриваются голоморфные функции в $D$.  Мы не используем  в доказательствах краткое сообщение \cite{KhaKha19} (в печати), дополняющее \cite{KhaRoz18},   и  некоторые специальные  частные версии результатов настоящей статьи из arXiv.org  \cite[\S~3]{KhT16}.

\subsection{Обозначения, определения, соглашения}\label{dfa} 
Всюду $\NN:=\{1,2, \dots  \}$ --- натуральные числа, $\RR\subset \CC$ --- {\it вещественная прямая,\/} $\RR^+:=\{x\in \RR\colon x\geq 0 \}$ --- {\it положительная полуось,\/} 
\begin{equation}\label{RR}
\RR^+_*:=\RR^+ \setminus \{0 \}, \quad 
\RR_{\pm \infty}:=\{ -\infty\} \cup \RR \cup \{+\infty \},  \quad 
\RR^+_{+\infty}:=\RR^+\cup \{+\infty \} 
\end{equation}  
с естественно дополненными неравенствами $-\infty\leq x \leq +\infty$ для любого $x\in \RR_{\pm\infty}$. 
Для    $r\in \RR^+_*$ и $z\in \CC$ полагаем $D(z,r):=\{z' \in \CC \colon |z'-z|<r\}$  --- {\it  открытый  круг с центром $z$ радиуса $r$;\/} $D(r):=D(0,r)$, $\DD:=D(1)$  и $D(z,+\infty):= \CC$. Для $z=\infty$ нам удобно принять $D(\infty,r):=\{z\in \CC_{\infty} \colon |z|>1/r\}$, $|\infty|:=+\infty$, и $D(\infty, +\infty):=\CC_{\infty}\setminus \{0\}$. Открытые круги $D(z,r)$ с $r\in \RR^+_*$ образуют открытую 
базу окрестностей точки $z\in  \CC_{\infty}$. Для $S\subset \CC_{\infty}$ через   $\Int S$, $\clos S$ и $\partial S$ обозначаем соответственно {\it внутренность, замыкание\/} и {\it границу\/} $S$  в $\CC_{\infty}$.  Для $S\subset S'\subset \CC_{\infty}$ пишем $S\Subset S'$, если  $ S$  --- относительно компактное подмножество  в   $S'$.  {\it (Под)область\/} в $\CC_{\infty}$ --- открытое связное подмножество в $\CC_{\infty}$. Всюду далее 
\begin{equation}\label{D}
D\neq \varnothing \text{\it  --- собственная подобласть в\/ $\CC_{\infty}\neq D$.}   
\end{equation}
Как и в \cite{KhaRoz18}, $\Har (S)$, $\sbh(S)$, $\dsbh(S)$,  $\Hol(S)$ и  $C^k(S)$ при $k\in \NN\cup\{ \infty\}$   --- классы соответственно  {\it гармонических, субгармонических\/} \cite{Rans}, \cite{L}, {\it $\delta$-субгармонических\/} \cite[3.1]{KhaRoz18}, {\it голоморфных,\/}  {\it $k$ раз непрерывно дифференцируемых\/} функций на  открытых множествах из $\CC_{\infty}$, включающих в себя $S\subset \CC_{\infty}$, но $C(S)$ --- класс {\it непрерывных\/} функций  именно на $S$.   
Через $\boldsymbol{-\infty}$ и $\boldsymbol{+\infty}$ обозначаем  функции, тождественно равные соответственно $-\infty$ и $+\infty$. В этих обозначениях  
\begin{equation}\label{dfF}
\sbh_*(S):=\sbh(S)\setminus \{\boldsymbol{-\infty}\}, \quad \sbh_*(S):=\dsbh(S)\setminus 
\{\boldsymbol{\pm\infty}\}, \quad \Hol_*(S):= \Hol(S)\setminus \{0\}.
\end{equation}
Символ $0$ --- нулевой вектор или начало отсчета в векторном или аффинном пространстве. Положительность в упорядоченном векторном пространстве $X$ всюду понимается как $\geq 0$; 
$+\infty\geq 0$ в\footnote{Метка-ссылка над знаками (не)равенства, включения, или, более общ\'о, бинарного отношения и т.\,п. означает, что данное соотношение как-то связано с отмеченной ссылкой.} $\RR_{+\infty}^+\overset{\eqref{RR}}{\subset} \RR^+_{\pm \infty}$. Для $A\subset X$  через $A^+$ обозначаем  множество положительных элементов из $A$. Класс всех функций $f\colon X\to Y$ обозначаем как $Y^X$.  Если $F(S)\subset \RR_{\pm \infty}^S:=(\RR_{\pm \infty})^S$ --- какой-либо класс {\it расширенных числовых функций,\/}  то $F^+(S)\subset (\RR^+_{+\infty})^S$ --- подкласс всех положительных функций из $F(S)$.

$\Meas  (S)$ --- класс\footnote{В \cite{KhaRoz18} использовалось обозначение $\mathcal M(S)$.} {\it борелевских вещественных мер\/} на борелевских подмножествах  множества  $S\subset  \CC_{\infty}$,       иначе называемых  {\it зарядами\/}  \cite{L}; $ \Meas _{\comp}(S)$  --- подкласс мер в $\Meas  (S)$ с компактным {\it носителем\/} $\supp \nu\Subset S$; $\Meas ^+ (S)$ --- положительные заряды, т.\,е. просто {\it меры;\/}
$\lambda$ --- {\it мера Лебега\/} в $\CC$,   $\delta_z$ --- {\it мера Дирака\/} в точке $z\in \CC_{\infty}$.

Пусть $ f\overset{\eqref{dfF}}{\in} \Hol_*(D)$. Функция $f$ {\it обращается в нуль на\/} последовательности точек\/ ${\mathsf Z}=\{{\mathsf z}_k\}_{k=1,2,\dots}$, лежащих в $D$ (пишем $\mathsf Z\subset D$), если кратность нуля, или корня,  функции $f$ в каждой точке $z \in D$ не меньше числа повторений этой точки  в последовательности ${\mathsf Z}$ (пишем $f({\mathsf Z})=0$).
Последовательности $\mathsf Z=\{{\mathsf z}_k\}_{k=1,2,\dots}\subset D$ {\it без предельных точек  в\/ $D$} сопоставляем
\begin{enumerate}[{\rm (1)}]
\item[{[div]}]\label{dz} {\it дивизор последовательности\/} $\sf Z$ на $D$ --- функция из $D$ в $\NN_0:=\{0\}\cup \NN$, равная в каждой точке $z\in D$ числу ее повторений в  $\sf Z$ и обозначаемая тем же символом  $\sf Z$, а именно:
\begin{equation}\label{divZ}
{\sf Z}(z):=\sum_{{\sf z}_k=z} 1=\sum_{k} \delta_{{\sf z}_k} (\{z \}), \quad z\in D;
\end{equation}
\item[{[cm]}]\label{nz}  {\it считающую   меру\/}
\begin{equation}\label{df:nZS}
n_{\mathsf Z}(S):=\sum_{{\mathsf z_k}\in S} 1=\sum_{k}\delta_{{\sf z}_k}(S),  \quad 
S\subset D, 
\end{equation}
--- число точек из $\mathsf Z$, попавших в $S$. Очевидно, ${\sf Z}(z)\overset{\eqref{divZ}}{\equiv}n_{\sf Z}(\{z\})$, $z\in D$. 
\end{enumerate}
Отходя от традиционной трактовки последовательности как функции натурального или целого аргумента, считаем две последовательности равными, если совпадают их дивизоры или, что эквивалентно, их считающие меры. Детальнее в   \cite[1.1]{Kh07}, \cite[0.1.2]{Khsur}.

{\it Последовательность нулей,\/} или корней, функции $f\in \Hol_* (D)$, каким-либо образом перенумерованную с учетом кратности, обозначаем через $\Zero_f $. Так как  $\ln  |f|\in \sbh (D)$ --- субгармоническая функция,  взаимосвязь ее меры Рисса $\nu_{\ln  |f|}$   при $f\neq 0$ со считающей мерой  ее нулей \eqref{df:nZS} задается равенством \cite[теорема 3.7.8]{Rans}, \cite[1.2.4]{KhaRoz18}
\begin{equation}\label{nufZ}
\nu_{\ln  |f|}=\frac1{2\pi} \Delta \ln  |f|\overset{\eqref{df:nZS}}{=}n_{\Zero_f}\in \Meas^+ (D), \quad \text{$\Delta$ --- {\it оператор Лапласа}.}
\end{equation}
В частности,  $f({\mathsf Z})=0$ эквивалентно  неравенству для мер $n_{\mathsf Z}\leq n_{\Zero_f}$ на $D$.

\subsection{Основные результаты}\label{111}
\begin{definition}[{\rm (версия понятия выметания \cite{Rans}, \cite{L}, \cite{BH})}]\label{dfbal}  Пусть 
 $S\Subset D$ и $F\overset{\eqref{RR}}{\subset} \RR_{\pm \infty}^{D\setminus S}$ --- некоторый класс расширенных числовых функций на $D\setminus S$.
Заряд $\mu\in \Meas(D)$ называем {\it аффинным выметанием\/} заряда $\nu\in \Meas( D)$ для $D$ вне $S\Subset D$ относительно $F$ и пишем $\nu\preceq_{S, F} \mu$, если найдется число $C\in \RR$, с которым 
\begin{equation}\label{bal}
\int_{D\setminus S}v\d \nu\leq \int_{D\setminus S} v \d \mu+C\quad \text{\it для всех $v\in F$},
\end{equation}
где интегралы в \eqref{bal}, вообще говоря, верхние \cite{B}. 
В частности, для последовательности ${\sf Z}=\{{\sf z}_k\}_{k=1,2,\dots}$ со считающей мерой $n_{\sf Z}$ из [cm]\eqref{df:nZS} заряд  $\mu\in \Meas(D)$ называем {\it аффинным выметанием\/} последовательности $\sf Z$ для $D$ вне $S\Subset D$ относительно класса $F $ и пишем ${\sf Z}\preceq_{S, F} \mu$, если $n_{\sf Z}\preceq_{S,F}\mu $, т.\,е. найдется число $C\in \RR$, с которым 
\begin{equation}\label{balZ}
\sum_{{\sf z}_k\in D\setminus S} v({\sf z}_k)\overset{\eqref{df:nZS}}{:=:}\int_{D\setminus S} v\d n_{\sf Z}\overset{\eqref{bal}}{\leq} \int_{D\setminus S} v \d \mu+C\quad \text{\it для всех $v\in F$}.
\end{equation}
\end{definition}
Очевидно, отношение {\it предпорядка\/} $\preceq_{S, F}$ на $\Meas(D)$ строго слабее стандартного отношения порядка $\nu\leq \mu$ на $\Meas (D)$.  Для функции $v\colon D\setminus S\to \RR_{\pm\infty}$ при $S\Subset D$ полагаем 
\begin{equation}\label{ldD}
\lim_{\partial D}v :=\lim_{D\ni z'\to z}v(z')\in \RR,\quad z\in \partial D,
\end{equation} 
если последний предел справа существует и один и тот же для любой точки $z\in \partial D$.

Во избежание некоторых чисто технических осложнений, связанных с необходимостью  применения  инверсии комплексной плоскости и преобразования Кельвина функций \cite[1.2.2]{KhaRoz18}, 
будем пока рассматривать только области $D\subset \CC$, т.\,е. $\infty \notin D$. 

\begin{theorem}[{\rm (критерий для субгармонических функций)}]\label{th1} Пусть 
$D\subset \CC$ --- область с неполярной границей $\partial D\subset \CC_{\infty}$,  
 $M\in \sbh_* (D)\cap C(D)$ --- функция с мерой Рисса $\mu\in \Meas^+(D)$,  
$\nu\in \Meas^+ (D)$, $b\in \RR^+_*$.
Тогда следующие три утверждения попарно эквивалентны. 
\begin{enumerate}[{\bf s1.}]
\item\label{s1}  Существует такая функция $u\in \sbh_*(D)$ с мерой Рисса $\nu_u\geq \nu$, что $u\leq M$ на $D$.  
\item\label{s2}  Для любого  подмножества $S\Subset D$, удовлетворяющего условиям 
\begin{equation}\label{SS}
\varnothing \neq \Int S\subset S=\clos S\Subset D, 
\end{equation}
 мера $\mu$ --- аффинное выметание меры $\nu$ для $D$ вне $S$ относительно класса тестовых субгармонических положительных функций\/ \cite[2.1]{KhaRoz18}
\begin{equation}\label{tv0}
\sbh^+_0(D\setminus S;\leq  b):=
\Bigl\{v\in \sbh^+(D\setminus S)\colon \lim_{\partial D}v \overset{\eqref{ldD}}{=}0, \; \sup_{D\setminus S}v\leq b  \Bigr\}.
\end{equation}
\item\label{s3}  Существует подмножество $S\Subset D$ из \eqref{SS}, для которого мера $\mu$ --- аффинное выметание меры $\nu$  для $D$ вне $S$ относительно класса 
$\sbh^+_{00}(D\setminus S;\leq  b)\cap C^{\infty}(D\setminus S)$, 
где 
\begin{equation}\label{tv00}
\sbh^+_{00}(D\setminus S;\leq  b):=
\bigl\{v\in \sbh_0^+(D\setminus S;\leq b)\colon \text{$\exists D_v\Subset D$,  
$v\equiv 0$ на $D\setminus D_v$} \bigr\}
\end{equation}
--- класс тестовых субгармонических финитных  функций\/ \cite[2.1]{KhaRoz18}, \cite[3.1]{KhT16}. 
\end{enumerate}
\end{theorem}
\begin{remark} Импликация {\bf s\ref{s2}}$\,\Rightarrow\,${\bf s\ref{s3}} 
очевидна для {\it любой\/} области $D\subset \CC_{\infty}$ и для {\it любой\/} меры
$\mu\in \Meas^+(D)$  без условия непрерывности $M$ на $D$.  
То же самое, как показано в  подразделе \ref{p12}, будет верно  и для импликации  {\bf s\ref{s1}}$\,\Rightarrow\,${\bf s\ref{s2}}. Лишь при доказательстве импликации {\bf s\ref{s3}}$\,\Rightarrow\,${\bf s\ref{s1}}  будет использованы как непрерывность $M$, так и неполярность  границы $\partial D\subset \CC_{\infty}$, что  эквивалентно существованию {\it функции Грина\/} $g_D$ для  области $D$ \cite[4.4]{Rans}, \cite[3.7, 5.7.4]{HK}. Теорема \ref{th:inver} из подраздела \ref{s31} --- утверждение  более  общее, чем импликация {\bf s\ref{s3}}$\,\Rightarrow\,${\bf s\ref{s1}}.  
\end{remark}

Перейдем к  голоморфной версии теоремы \ref{th1}, в которой возникает некоторый зазор между необходимыми и достаточными условиями. Этот, хотя и незначительный, зазор вряд ли может быть ликвидирован даже для круга  $D=\DD$ в рассматриваемой здесь  общей ситуации. 
Через $\dist (\cdot,\cdot )$ обозначаем {\it функцию евклидова расстояния\/} между парами точек, между  точкой и множеством, между множествами в $\CC$. По определению полагаем $\dist (\cdot , \varnothing):=:\dist (\varnothing,\cdot ):=\inf \varnothing:=+\infty =:\dist (z,\infty):=:\dist (\infty, z)$ при  $z\in \CC$.

\subsubsection{Выбор поднятия функции $M$}\label{ssM}
Далее $r\colon D\to (0,1]$ --- произвольная непрерывная функция, удовлетворяющая условию  
\begin{equation}\label{rdD}
0<r(z)<\dist(z, \partial D) \quad\text{\it при  всех $z\in D$.}
\end{equation}  
$L^1_{\loc} (D)\subset \RR_{\pm \infty}^D$ --- класс {\it локально интегрируемых по мере Лебега\/} $\lambda $ функций. Функции $M\in L^1_{\loc} (D)$ сопоставляем ее переменные усреднения  по  кругам 
\begin{multline}\label{M*r}
M^{*r}  (z)
:=\frac{1}{\lambda(D(z,r(z)))}\int_{D(z,r(z))} M \d \lambda\\
=\frac{1}{\pi r^2(z)}\int_0^{2\pi}\int_0^{r(z)} M(z+te^{i\theta})\, t\d t \d \theta, \quad D(z,r(z))\subset D,
\end{multline}
 а также, следуя \cite{Kha15}, \cite{KhaKha16},    ее <<{\it поднятие\/}>> 
$M^{\uparrow}$, а именно: 
\begin{enumerate}[{(i)}]
\item\label{Di}  в общем случае $D\subset \CC$ полагаем
\begin{equation}\label{Dg}
M^{\uparrow}(z):=M^{*r}(z)+\ln\frac{1}{r(z)} +(1+\varepsilon )\ln (1+|z|)\quad\text{для всех  $z\in D$},
\end{equation}
где число  $\varepsilon\in \RR^+_*$ может быть выбрано сколь угодно малым;
\item\label{Dii} если $\CC_{\infty}\setminus \clos D\neq \varnothing$ или область $D\subset \CC$ односвязная в $\CC_{\infty}$, то  
\begin{equation}\label{df:Nd}
M^{\uparrow}(z):= M^{*r}(z)+\ln\frac{1}{r(z)} \quad\text{для всех  $z\in D$};
\end{equation}
\item\label{Diii} если $D=\CC$, то для любого сколь угодно большого числа  $P>0$ можем положить
\begin{equation}\label{df:Nd_C}
M^{\uparrow}(z):= M^{*r}(z), \quad  r(z):=\frac{1}{(1+|z|)^P}\quad\text{для всех  $z\in D$}.
\end{equation}
\end{enumerate}

\begin{theorem}[{\rm (необходимые/достаточные условия для голоморфных функций)}]\label{th2}
Пусть $D$ --- область в $\CC$, $M\in \sbh_*(D)$ с мерой Рисса $\mu\in \Meas^+(D)$,  ${\sf Z}=\{{\sf z}_k\}_{k=1,2,\dots}\subset D$, $b\in \RR^+_*$. Каждое из следующих трех утверждений 
{\bf h\ref{h1}}--{\bf h\ref{h3}} следует из предыдущего. 
\begin{enumerate}[{\bf h1.}]
\item\label{h1}  Существует функция  $f\overset{\eqref{dfF}}{\in} \Hol_*(D)$, для которой  $f({\sf Z})=0$ и    $|f|\leq \exp M$ на $D$.  
\item\label{h2}  Для любого множества $S$ из \eqref{SS} мера $\mu$ --- аффинное выметание последовательности  $\sf Z$ для $D$ вне $S$ относительно класса тестовых субгармонических функций\/ \eqref{tv0}.

\item\label{h3}  Существует множество $S$ из \eqref{SS}, для которого $\mu$ --- аффинное выметание последовательности  $\sf Z$  для $D$ вне $S$ относительно класса
$\sbh^+_{00}(D\setminus S;\leq  b)\cap C^{\infty}(D\setminus S)$.
\end{enumerate}
Обратно, когда дополнительно  $\partial D$ --- неполярное множество в $\CC_{\infty}$ и 
$M\in C(D)$, утверждение  {\bf h\ref{h3}} влечет за собой существование функции $f\in \Hol_*(D)$, обращающейся в нуль на $\sf Z$ и удовлетворяющей неравенству $|f|\overset{\ref{ssM}}{\leq} \exp M^{\uparrow}$ на $D$ для $M^{\uparrow}$ из\/ \eqref{Di}--\eqref{Diii}. 
\end{theorem}
 \begin{remark}\label{r2}
Импликация {\bf h\ref{h2}}$\,\Rightarrow\,${\bf h\ref{h3}} очевидна.  В \S\S~\ref{Grf},\,\ref{apnd} даются теоремы обращения \ref{th3} и \ref{th:4} исключительно в терминах аффинного выметания относительно классов соответственно  функций Грина и определенных логарифмических потенциалов аналитических дисков. Следствия \ref{holDg} и \ref{holDgad}   из теорем \ref{th3} и \ref{th:4}   дают иные варианты импликации {\bf h\ref{h3}}$\,\Rightarrow\,${\bf h\ref{h1}}. 
\end{remark}

\section{Доказательства теорем \ref{th1} и \ref{th2}}\label{ss2}

\subsection{Доказательство импликаций {\bf s\ref{s1}}$\,\Rightarrow\,${\bf s\ref{s2} и 
{\bf h\ref{h1}}$\,\Rightarrow\,${\bf h\ref{h2}}}}\label{p12}
В этом подразделе не предполагается, что функция $M$ непрерывна. Непустая область $D\subset \CC$  произвольная. 

Выберем $z_0\in D$ так,что $u(z_0)\neq -\infty$ и $M(z_0)\neq -\infty$.  Выбор регулярной области $\widetilde{D}$, $S\Subset \widetilde D \Subset D$, участвующей в формулировке \cite[основная теорема]{KhaRoz18}, произволен. Тогда из условия   {\bf s\ref{s1}} по \cite[основная теорема]{KhaRoz18} найдутся числа $C,\overline C_M\in \RR^+$, для которых \cite[(3.3)]{KhaRoz18}
\begin{equation*}
C u(x_0) 	+\int_{D\setminus S}  v \d {\nu}_u \leq	\int_{D\setminus S}  v \d \mu	 +C\, \overline{C}_M
\quad\text{для всех $v\in \sbh^+_0(D\setminus S;\leq  b)$}.
\end{equation*} 
Это  ввиду $\nu \leq \nu_u$ по определению \ref{dfbal} показывает, что 
$\nu \preceq_{S,F} \nu_u \preceq_{S,F} \mu$  для аффинного выметания $\preceq_{S,F}$ для $D$ вне $S$ относительно  класса $F\overset{\eqref{tv0}}{=}\sbh^+_0(D\setminus S;\leq  b)$. 

Импликация {\bf h\ref{h1}}$\,\Rightarrow\,${\bf h\ref{h2}} --- частный случай 
импликации {\bf s\ref{s1}}$\,\Rightarrow\,${\bf s\ref{s2}} для $u:=\ln|f|$ с мерой Рисса 
$n_{\Zero_f}\overset{\eqref{nufZ}}{\geq} n_{\sf Z}$ в рамках определения \ref{balZ} аффинного выметания $\mu$ последовательности ${\sf Z}$  для $D$ вне $S$ через  соотношение \eqref{balZ}.

\subsection{Доказательство импликации {\bf s\ref{s3}}$\,\Rightarrow\,${\bf s\ref{s1}}}\label{s31}

Будет доказано более общее утверждение.

\begin{theorem}\label{th:inver}  
Пусть $D\subset \CC_{\infty}$ --- область с неполярной границей $\partial D$ и $S\Subset D$ --- подмножество из \eqref{SS},  $M\in \dsbh_* (D)$ --- $\delta$-субгармоническая  функция с зарядом Рисса $\nu_M\in \Meas (D)$,  $\nu\in \Meas^+(D)$,  $b\in \RR_*^+$. Если заряд $\nu_M$  --- аффинное выметание меры  $\nu$ для $D$ вне $S$ относительно класса  $\sbh_{00}^+(D\setminus S;\leq b)\cap C^{\infty}(D\setminus S)$, т.\,е. существует число  $C\in \RR$, с которым 
\begin{equation}\label{estv+0}
\int\limits_{D\setminus S}  v \dd {\nu} 		\leq	\int\limits_{D\setminus S}  v \dd {\nu}_M	+C
\quad \text{для всех $v\overset{\eqref{tv00}}{\in}  \sbh_{00}^+(D\setminus S;\leq b)\cap C^{\infty}(D\setminus S)$}, 
\end{equation}
то  для каждой непрерывной функции\/  $r \colon D\to (0,1]$, удовлетворяющей условию\/ 
\eqref{rdD},  найдется функция  $u\in \sbh_*(D)$ с мерой Рисса $\nu_u\geq \nu$, для которой
\begin{equation}\label{in:hM}
u(z)\overset{\eqref{M*r}}{\leq}  M^{*r} (z) \quad\text{для всех  $ z\in D$}. 
\end{equation}
Если, в дополнение, $M\in \dsbh_*(D)\cap C(D)$, т.\,е. $M$ еще и непрерывна, то 
при условии \eqref{estv+0} можно подобрать функцию  $u\in \sbh_*(D)$ с мерой Рисса $\nu_u\geq \nu$   так, что   $u\leq M$  на $D$.
\end{theorem}
\begin{proof}  Будем сначала предполагать, что вместо условия  \eqref{estv+0} для некоторой непустой подобласти 
\begin{equation}\label{D0S}
D_0\Subset  \Int S\Subset S\overset{\eqref{SS}}{=}\clos S \Subset D,
\end{equation} 
и некоторого числа $C\in \RR$ имеет место неравенство 
\begin{equation}\label{estv+}
\int\limits_{D\setminus D_0}  v \dd {\nu} 		\leq	\int\limits_{D\setminus D_0}  v \dd {\nu}_M	+C
\quad \text{для всех $v\overset{\eqref{tv00}}{\in}  \text{\rm sbh}_{00}^+(D\setminus D_0;\leq b)$}, 
\end{equation}
т.\,е. тестовые финитные функции $v$ не обязательно дифференцируемы, а $S\Subset D$ несколько сужено  до подобласти  $D_0\overset{\eqref{estv+}}{\Subset} D$. Для  меры $\nu$ на $D\supset D_0$ всегда можно  подобрать точку $z_0\in  D_0$ так, что корректно определено значение $M(z_0)\neq \pm\infty $, т.\,е. $z_0\in D_0\cap \dom M$  в обозначении из \cite[3.1]{KhaRoz18}, а также  для некоторого числа $r_0>0$ выполнено 
 \begin{equation}\label{dMa+}
\left(\int_0^{r_0}\frac{\nu (z_0,t)}{t} \dd t<+\infty\right)\; {\Longleftrightarrow}	\;
\left(\int_{D(z_0,r_0)} \ln |z'-z_0| \dd \nu (z')>-\infty\right), \quad D(z_0,3r_0)\Subset D_0.
\end{equation}
Условия \eqref{dMa+}, в частности, обеспечивают существование    функции $u_0\in \sbh_*(D)$ с мерой Рисса  $\nu_{u_0}=\nu$  и свойством $u_0(z_0)\neq -\infty$ \cite[3.1]{KhaRoz18}. 

Далее нам временно потребуется ограниченность функции $M$ в окрестности точки $z_0$. Для этого пока преобразуем  ее локально с сохранением условия 
\eqref{estv+}. Из \eqref{dMa+} и представления  $M=u_+-u_-$ в  виде разности субгармонических функций $u_+, u_- \in \sbh_*(D)$ можно локально изменить значения функции $M$ в  $D(z_0,2r_0)\Subset D_0$, а именно: гармонически  продолжить интегралом Пуассона    функции $u_+$ и $u_-$ внутрь  $D(z_0,2r_0)$. Обозначаем их соответственно $u_+^{\circ}$ и $u_-^{\circ}$. Тогда 	$M^{\circ}:=u_+^{\circ}-u_- ^{\circ}\in \dsbh_* (D)$ --- ограниченная функция в окрестности замкнутого круга 	$\clos D(z_0,r_0)$, а \eqref{estv+} по-прежнему выполнено для   всех  $v{\in}  \text{\rm sbh}_{00}^+(D\setminus D_0;\leq b)$. Пока будем обозначать   функцию $M^{\circ}$ прежним символом $M$. 
Через $J_{z_0}(D)$ обозначаем, как и в \cite{KhaRoz18}, класс всех мер Йенсена $\mu \in \Meas_{\comp}^+(D)$, удовлетворяющих условию $u(z_0)\leq \int u\d \mu$ для всех $u\in \sbh(D)$.
 Будет использована
\begin{theoA}[{\rm (частный случай \cite[теорема 6]{Kh07}})]
Пусть\/ $M\in L^1_{\loc} (D)$,  $z_0\in D$,   $u_0 \in \sbh(D)$    с $u_0(z_0)\neq -\infty$.
Если  функция $M$ ограничена в открытой окрестности замыкания\/ 
$\clos D_1$ какой-нибудь подобласти\/ $D_1\Subset D$, содержащей $z_0$, и  существует число\/ $C_0\in \RR$, с которой  
\begin{equation}\label{in:arsV}
	\int_{D} u_0 \dd \mu \leq \int_{D} M \dd \mu +C_0 
\quad\text{для любой  меры\/ $ \mu \in  J_{z_0}(D)$,}	
\end{equation}
 то  для каждой непрерывной функции\/  $r \colon D\to \RR_*^+$, удовлетворяющей условию\/ \eqref{rdD},
	 найдется функция\/ $w\in \sbh_*(D)$, для которой
\begin{equation}\label{in:hMsh}
u_0+w\leq  M^{*r} \quad\text{на  $ D$}.
\end{equation}
\end{theoA}
В нашем случае роль области $D_1$ будет играть круг $D(z_0,r_0)$. Кроме того, потребуются

\subsubsection{Потенциалы Йенсена}\label{mpj}
 Функцию $V\in \sbh^+\bigl(\CC_{\infty} \setminus \{z_0\}\bigr)$ 
называем {\it потенциалом Йенсена  внутри  $D$ с полюсом в  $z_0\in D$} \cite[определение 3]{KhaRoz18}, если выполнены два условия:
\begin{enumerate}[{\rm 1)}]
\item\label{V:f2} {\it найдется  область $D_V\Subset D$, содержащая  $z_0\in D_V$, для которой $V (z)\equiv 0$ при   $z\in \CC_{\infty} \setminus D_V$;  
\item\label{V:f3}  имеет место  {\it логарифмическая полунормировка в точке\/ $z_0$,} а именно: 
\begin{subequations}\label{nvz}
\begin{align}
\limsup\limits_{z_0\neq z \to z_0}\dfrac{V(z )}{l_{z_0}(z)}&\leq  1, 
 \tag{\ref{nvz}o}\label{nvzn}\\ 
\text{\it где}\quad l_{z_0}(z):=&\begin{cases}
\ln \frac1{|z-z_0|} \quad&\text{\it при $z_0\neq \infty$},\\
\ln |z| \quad&\text{\it при $z_0= \infty$}.
\end{cases}
\tag{\ref{nvz}l}\label{{nvz}INF}
\end{align}
\end{subequations}}
\end{enumerate}
Класс всех  таких потенциалов Йенсена обозначаем через  $PJ_{z_0} (D)$.  

{\it Логарифмический потенциал рода\/ $0$ вероятностной меры\/ $\mu \in \mathcal  M^+_{\comp}(\CC_{\infty})$  с полюсом в точке $z_0\in \CC_{\infty}$} определяем для всех $w\in \CC_{\infty}\setminus \{z_0\}$ как функцию
\begin{subequations}\label{df:Vmu}
\begin{align}
V_{\mu}(w ) 	&:= 		\int_{D} \ln \Bigl|\frac{w-z}{w-z_0}\Bigr| \dd \mu (z)
= \int_{D} \ln \Bigl|1-\frac{z-z_0}{w-z_0}\Bigr| \dd \mu (z) \quad\text{\it при $z_0\neq \infty$},
\tag{\ref{df:Vmu}o}\label{{df:Vmu}o}
\\
\intertext{где при $w=\infty$ подынтегральные выражения доопределены значением  $0$,} 
V_{\mu}(w ) 	&:=	\int_{D} \ln \Bigl|\frac{w-z}{z}\Bigr| \dd \mu (z)=
\int_{D} \ln \Bigl|1-\frac{w}{z}\Bigr| \dd \mu (z) \quad\text{\it при $z_0= \infty$},
\tag{\ref{df:Vmu}$\infty$}\label{{df:Vmu}in}
\end{align}
\end{subequations}
где при $z=\infty$ подынтегральные выражения доопределены значением  $0$.

Напомним основные взаимосвязи между $J_{z_0} (D)$ и $PJ_{z_0} (D)$. 
Первая --- следующее утверждение о двойственности.
\begin{propos}[{\cite[предложение 1.4, теорема двойственности]{Khab03}}]\label{pr:1} Отображение  
\begin{equation*}
	\mathcal P \colon J_{z_0}(D)\to PJ_{z_0} (D), \quad 
	\mathcal P (\mu)\overset{\eqref{df:Vmu}}{:=}  V_{\mu}, \quad \mu \in  J_{z_0}(D), 
\end{equation*}
---  биекция, для которой 
$\mathcal P\bigl(t\mu_1+(1-t)\mu_2\bigr)=t\mathcal P (\mu_1)+(1-t)\mathcal P (\mu_2)$ для всех $t\in [0,1]$, 
 а обратное отображение ${{\mathcal P}}^{-1}$ определено равенством
\begin{equation}\label{eq:mu}
	{{\mathcal P}}^{-1}(V)\overset{\eqref{{nvz}INF}}{=}\frac1{2\pi} \Delta  V\Bigm|_{D\setminus \{z_0\}}+
	\Bigl(1-\limsup\limits_{z_0\neq z \to z_0}\dfrac{V(z )}{l_{z_0}(z)}\Bigr)\cdot
	{\delta}_{z_0}\, , \quad V\in PJ_{z_0}(D).
\end{equation}
\end{propos}

Вторая   ---  это {\it расширенная формула Пуассона--Йенсена\/} \eqref{f:PJ}. 
\begin{propos}[{\cite[предложение 1.2]{Khab03}}]\label{pr:2}
Пусть $\mu \in J_{z_0}(D)$. Тогда для любой функции $u\in \sbh (D)$ с мерой  Рисса
$\nu_u$ при $u(z_0)\neq -\infty$ имеем равенство
\begin{equation}\label{f:PJ}
	u(z_0) +\int_{D\setminus \{z_0\}} V_{\mu} \dd {\nu}_u=\int_{D} u \dd \mu  .
\end{equation}
\end{propos}

\begin{lemma}\label{pr:PJP} Пусть $M\in \dsbh_*(D)$ с зарядом Рисса $\nu_M$,  $z_0\in\dom M$,   $u_0\in \sbh (D)$ с мерой Рисса $\nu$,
 $u_0(z_0)\neq -\infty$,  а $V\in PJ_{z_0}(D)$ --- потенциал Йенсена  и $C_1\in \RR$. Если 
\begin{equation}\label{inmuM+}
	\int_{D\setminus \{z_0\}} V\dd\nu \leq \int_{D\setminus \{z_0\}} V \dd\nu_M +C_1, 
\end{equation}
то для меры Йенсена  $\mu\overset{\eqref{eq:mu}}{=}\mathcal P^{-1} (V)\in J_{z_0}(D)$ справедливо неравенство
\begin{equation}\label{inmuM+q}
	\int u_0\dd\mu \leq \int M \dd\mu +C_0, \quad \text{где $C_0=C_1-M(z_0)+u_0(z_0)$}.
\end{equation}
\end{lemma}
\begin{proof}[леммы \ref{pr:PJP}]
При условии $z_0\in \dom M$ функция $M$ представима разностью  $M=u_+-u_-$  функций $u_{\pm} \in \sbh_*(D)$  с мерами Рисса соответственно $\nu_M^{\pm} \in \Meas^+(D)$, для которых  имеем $u_{\pm}(z_0)\neq -\infty$. К каждой из функций $u_{\pm} $ применима расширенная  формула Пуассона\,--\,Йенсена предложения \ref{pr:2}, следовательно, она применима и к функции $M$. Отсюда  для  меры Йенсена  $\mu\overset{\eqref{eq:mu}}{:=}\mathcal P^{-1}(V)$ имеем
 \begin{multline*}
	\int_{D} u_0 \dd \mu   \overset{\eqref{f:PJ}}{=}\int_{D\setminus \{z_0\}} V\dd {\nu} +u_0(z_0) \\
	\overset{\eqref{inmuM+}}{\leq} \int_{D\setminus \{z_0\}} V \dd\nu_M +C_1+u_0(z_0)
	\overset{\eqref{f:PJ}}{=} \int M \dd\mu -M(z_0)+C_1 +u_0(z_0),
\end{multline*}
что и доказывает требуемое неравенство  \eqref{inmuM+q}.
\end{proof}

Вернемся непосредственно к доказательству теоремы \ref{th:inver}.  
Для области $D$ с неполярной границей  $\partial D\subset \CC_{\infty}$ всегда существует функция Грина $g_{D}(\cdot ,z_0)$ 
{\; \it  с полюсом в точке $z_0$.} Далее всюду  в нашем доказательстве для краткости записей 
\begin{equation*}
g:=g_D(\cdot, z_0) \quad \text{\it --- функция Грина для $D$ с полюсом $z_0\in D_0$.}
\end{equation*}
Здесь для нас важны только  следующие ее свойства {\large(}см. \cite[4.4]{Rans}, \cite[3.7, 5.7]{HK}{\large)}:
\begin{enumerate}[{(g1)}]
	\item\label{g1} $\lim\limits_{z_0\neq z\to z_0} \dfrac{g(z )}{l_{z_0}(z)}\overset{\eqref{{nvz}INF}}{=}  1$ --- нормировка в точке $z_0$ \eqref{nvzn};
	\item\label{g2} $g\in \Har^+\bigl(D\setminus \{z_0\}\bigr)$ --- гармоничность и положительность в $D\setminus \{z_0\}$. 
	\end{enumerate}
В частности, из принципа максимума-минимума,  ввиду $z_0\in D_0\Subset D$, 
\begin{equation}\label{B0}
	0<\const_{z_0,D_0,D}:=B_0:=\sup_{z\in \partial D_0} g(z)<+\infty.
\end{equation}
Пусть $V\in PJ_{z_0}(D)$ --- {\it произвольный потенциал Йенсена.\/} Тогда ввиду (g\ref{g1})--(g\ref{g2}) 
  	\begin{equation*}
	\limsup_{D\ni z\to z_0} \frac{(V-g)(z)}{l_{z_0}(z)} \overset{\rm(g\ref{g1})}{\leq} 0,  \quad 
	V-g \overset{\rm(g\ref{g2})}{\in} {\sbh_*} \bigl(D\setminus \{z_0\}\bigr).	
	\end{equation*}
Отсюда точка $z_0$ --- устранимая особенность для функции $V-g\in \sbh_*(D\setminus \{z_0\})$ и, поскольку
\begin{equation*}
\limsup_{D\ni z'\to z} (V-g)(z')\leq \limsup_{D\ni z'\to z} V(z')  = 0 \quad	\text{для всех $z\in \partial D$,}
\end{equation*}
для функции $V-g\in \sbh_*(D)$ по принципу максимума $V-g\leq 0$ на $D$, т.\,е. 
\begin{equation}\label{es:g}
V\leq g \quad \text{на $D$}, \qquad 
V\overset{\eqref{B0}}{\leq} B_0 \quad \text{на $\partial D_0$}. 	
\end{equation}
Следовательно, для рассматриваемой в открытой окрестности $D\setminus D_0$  функции 
\begin{equation*}
v:=\frac{b}{B_0} \, V	\in \sbh_{00}^+\bigl(D\setminus D_0;\leq b\bigr)
\end{equation*}
справедливо неравенство \eqref{estv+}. Умножая обе его части на $B_0/b$,  получаем 
\begin{equation*}
\int_{D\setminus D_0} V \dd \nu  	\leq \int_{D\setminus D_0} V \dd\nu_M +\frac{B_0}{b} \,C 
\quad\text{\it для всех $V\in PJ_{z_0}(D)$.}
\end{equation*}
Это неравенство можно переписать в виде
\begin{multline}\label{bQ}
\int_{D\setminus \{z_0\}} V \dd \nu  \leq \int_{D\setminus \{z_0\}} V \dd\nu_M 
+\frac{B_0}{b} \,C  + \int_{D_0\setminus \{z_0\}} V \dd \nu  +\int_{D_0\setminus \{z_0\}} V \dd\nu^-_M \\
\overset{\eqref{es:g}}{\leq}  \int_{D\setminus \{z_0\}} V \dd\nu_M 
+\left(\frac{B_0}{b} \,C  + \int_{D_0\setminus \{z_0\}} \bigl(g \dd \nu  +
g \dd\nu^-_M\bigr)\right) 	\quad\text{\it для всех $V\in PJ_{z_0}(D)$.}
\end{multline}
Последний парный интеграл здесь конечен ввиду \eqref{dMa+} и 	$z_0\in D_0\cap \dom M$, а также не зависит от $V\in PJ_{z_0}(D)$. Таким образом, с постоянной $C_1$, равной  значению <<большой>> скобки в правой части \eqref{bQ},  выполнено  \eqref{inmuM+}	 для любого потенциала $V\in PJ_{z_0}(D)$. Отсюда по лемме  \ref{pr:PJP} имеет место \eqref{inmuM+q} для любой меры Йенсена $\mu\in J_{z_0}(D)$. Следовательно, выполнено  условие   \eqref{in:arsV}  теоремы A и найдется  функция  $w\in \sbh_* (D)$, для которой имеем \eqref{in:hMsh}. При этом c мерой Рисса $\nu_w$ функции $w$, очевидно, выполнено  неравенство $\nu+\nu_w\geq \nu$ на $D$. Следовательно, функция $u^{\circ}:=u_0+w$ с мерой Рисса $\nu_{u^{\circ}}=\nu+\nu_w$ --- требуемая в \eqref{in:hM}, но пока для функции $M=M^{\circ}$, отличающейся от $M$ в круге 
$D(z_0,2r_0)$. Вернемся к прежним обозначениям $M\leq M^{\circ}$.   Для {\it непрерывной\/} функции  $r$ функции  $M^{*r}$, $(M^{\circ})^{*r}$ также непрерывны в $D$, поскольку обе они из  класса  $L_{\loc}^1(D)$. В то же время субгармоническая функция $u^{\circ}\neq \boldsymbol{-\infty}$ ограничена сверху в $D(z_0,3r_0)\Subset D$. Следовательно,  можно выбрать достаточно большую  постоянную $C_2\geq 0$ так, что $u_0:=u^{\circ}-C_2\leqslant (M^{\circ})^{*r}$ на $D$  с мерой Рисса $\nu_{u_0}=\nu_{u^{\circ}}\geq \nu$. По условиям  на функцию $r$  найдется подобласть $D_2\Subset D$, включающая в себя $D(z_0,r_0)$,  для которой по построению $M^{\circ}$ и определению усреднения
на $D\setminus D_2$ выполнено равенство $(M^{\circ})^{*r}=M^{*r}$, а значит и неравенство  $u_0\leqslant M^{*r}$ на $D\setminus D_2$. Поскольку функция $M^{*r}$ непрерывна на $D$, а $u_0$ ограничена сверху на $D_1$, можно выбрать достаточно большое  число $C_3\geq 0$ так, что $u:=u_0-C_3\leq M^{*r}$ на $D$ с мерой Рисса $\nu_u=\nu_{u_0}\geq \nu$, что и дает  \eqref{in:hM}.

Если функция $M\in \dsbh_*(D)$ непрерывна, то она  изначально локально ограничена снизу, что позволяет избежать в доказательстве промежуточного использования функции $M^{\circ}$. Кроме того, непрерывная функция $M$ локально равномерно непрерывна, что позволяет  выбрать непрерывную функцию $r$, удовлетворяющую условию \eqref{rdD}, с которой $M^{*r}\leq M+1$ на $D$. Это дает возможность  заменить  правую часть $M^{*r}$ в \eqref{in:hM} на $M$.

Пусть теперь выполнено условие  \eqref{estv+0} доказываемой теоремы \ref{th:inver}.  Выведем из него  условие  \eqref{estv+}, при котором все заключения теоремы \ref{th:inver} уже доказаны. 

 Пусть   $v \in  \text{\rm sbh}_{00}^+(D\setminus D_0;\leq b)$ --- {\it произвольная финитная\/} тестовая функция, где подобласть $D_0\Subset \Int S$ выбрана как в \eqref{D0S}.  Ввиду финитности функции $v$ найдется подобласть   $D_v \Subset D$, для которой  $D_0\Subset D_v$, а  функция  $v$ субгармоническая на $D\setminus D_0$ и тождественно равна нулю на $D\setminus D_v$. При этом  положим 
\begin{equation}\label{e0}
\e_0:=\frac{1}{2}\min \bigl\{ \dist (D_v, \partial D ), 
\; \dist (D_0, D\setminus S) \bigr\}>0.
\end{equation}
Рассмотрим  бесконечно дифференцируемую функцию $a\colon \RR^{+} \to \RR^+$ с {\it носителем}   
\begin{equation}\label{a}
\supp a\subset (0,1) \quad\text{и {\it нормировкой}}\quad  2\pi \int_0^{+\infty} a(x)x\dd x =1,
\end{equation} 
а также  меры $\alpha_{\e}\in \Meas^+(\CC)$, определяемые плотностями
\begin{equation}\label{ae}
	\dd \alpha_{\e} (z)\overset{\eqref{a}}{:=}\frac{1}{\e^2}\,a\bigl(|z|/\e\bigr)\dd \lambda (z), \quad 0<\e< \e_0, \quad z\in \CC.
\end{equation}
Как известно \cite[2.7]{Rans}, \cite[3.4.1]{HK},  для убывающей последовательности строго положительных чисел $\e_n\underset{n\to \infty}{\longrightarrow} 0$, $\e_n\overset{\eqref{e0}}{\leq} \e_0$,  
последовательность субгармонических бесконечно дифференцируемых функций-сверток  
$	v_n:=v*\alpha_{\e_n}$, убывая по $n\in \NN$, поточечно стремиться к функции $v$ на $D\setminus S$. В частности, согласно \eqref{e0},  $v\leq v_n$ на $D\setminus S$ и $v_n\leq b$ на $D\setminus S$ как усреднения по вероятностным мерам \eqref{ae}.  По построению все функции $v_n\in  {\sbh}_{00}^+(D\setminus S;\leq b)\cap C^{\infty} (D\setminus S)$.   По условию  \eqref{estv+0} найдется число $C$, с которым 
\begin{equation*}
\int\limits_{D\setminus S}  v_n \dd {\nu} 		\leq	\int\limits_{D\setminus S}  v_n \dd {\nu}_M	+C
\end{equation*}	
  для всех построенных функций $v_n$ при всех  $n\in \NN$. 
	Отсюда по  разложению Хана\,--\,Жордана	для заряда Рисса $\nu_M	=\nu_M^+-\nu_M^-$, $\nu_M^{\pm}\in \Meas^+(D) $,  имеем
	\begin{equation*}
\int\limits_{D\setminus S}  v_n \dd ({\nu}+\nu_M^-) 		\leq	\int\limits_{D\setminus S}  v_n \dd {\nu}_M^+	+C'
\quad\text{для всех $n\in \NN$},
\end{equation*}
что  ввиду $v\leq v_n$ на $D\setminus S$ дает 
\begin{equation*}
\int\limits_{D\setminus S}  v \dd {\nu} 		\leq	\int\limits_{D\setminus S}  v_n \dd {\nu}_M^+ - 
\int\limits_{D\setminus S}  v \dd \nu_M^-	+C\quad\text{для всех $n\in \NN$}.
\end{equation*}	
Устремляя  в первом интеграле справа  $n$ к $+\infty$, 	 ввиду убывания последовательности тестовых финитных бесконечно дифференцируемых  функций  $v_n$ к  $v\in \sbh_{00}^+(D\setminus D_0;\leq b)$ имеем 
\begin{equation}\label{prom}
\int\limits_{D\setminus S}  v \dd {\nu} 		\leq	\int\limits_{D\setminus S}  v \dd {\nu}_M^+ - 
\int\limits_{D\setminus S}  v \dd \nu_M^-	+C=
\int\limits_{D\setminus S}  v \dd {\nu}_M+C.
\end{equation}	
Определим постоянные $C_4,C_5\in \RR^+$, не зависящие  от $v$, соотношениями  
\begin{equation*}
0\leq \int_{D_0\setminus S}v\d \nu\leq b\,\nu(D_0\setminus S)=:C_4<+\infty,\quad
0\leq \int_{D\setminus S} v \d \nu_M^-\leq b\,\nu_M^-(D_0\setminus S)=:C_5<+\infty
\end{equation*}
Тогда неравенство  \eqref{prom} без промежуточной разности интегралов останется в силе, если 
интегрирования по $D\setminus S$ заменить на интегрирования по $D\setminus D_0$ с заменой постоянной $C$  на постоянную $C+C_4+C_5$
 В силу произвола в выборе  тестовой финитной функции  $v{\in}  \text{\rm sbh}_{00}^+(D\setminus D_0;\leq b)$ неравенство \eqref{estv+}  с новой постоянной $C+C_4+C_5$ вместо $C$ выполнено для всех таких $v$. Это  завершает 
доказательство теоремы \ref{th:inver}. 
\end{proof}
\begin{remark}\label{add:2}  В случае  функции   $M\in\sbh_*(D)$  в теореме\/ {\rm \ref{th:inver}} достаточно требовать, чтобы  функция $r$ со свойством    \eqref{rdD}   была лишь локально отделенной от нуля снизу в  том смысле, что для любого $z\in D$ найдется число $t_z >0$, для которого
$D(z,t_z)\Subset D$ и $\sup\limits_{D(z,t_z)}r >0$. Действительно, из элементарных геометрических соображений, использующих компактность, например, с использованием исчерпания области $D$ последовательностью относительно компактных подобластей, устанавливается
\begin{lemma}\label{l:dz}
Для отделенной от нуля снизу функции $r$ на $D$, удовлетворяющей условию\/ \eqref{rdD}, найдется непрерывная и даже бесконечно дифференцируемая функция $\hat r\leq r$, по-прежнему удовлетворяющая условию\/ \eqref{rdD}.
\end{lemma}
Применяя теорему \ref{th:inver} с непрерывной  функцией $\hat r$ вместо $r$, строим необходимую   функцию $u\overset{\eqref{in:hM}}{\leq} M^{*{\hat r}}\leq M^{*r}$, где  использовано возрастание  усреднений \eqref{M*r} по $r$ для  $M\in \sbh_*(D)$. 
\end{remark}

\subsection{Доказательство импликации {\bf h\ref{h3}}$\,\Rightarrow\,${\bf h\ref{h1}} с поднятием $M^{\uparrow}$}\label{h31}

В условиях части <<{\it Обратно, \dots}>> теоремы \ref{th2}  из утверждения  {\bf h\ref{h3}} и импликации {\bf s\ref{s3}}$\,\Rightarrow\,${\bf s\ref{s1}} теоремы \ref{th1} следует  существование субгармонической функции $u\in \sbh_*(D)$ с мерой Рисса $\nu_u\geq n_{\sf Z}$, удовлетворяющей неравенству $u\leq M$ на $D$. Поскольку по теореме Вейерштрасса всегда существует голоморфная функция $f_{\sf Z}\in \Hol_*(D)$ с последовательностью нулей $\Zero_{f_{\sf Z}}={\sf Z}$, то последнее означает, что существует функция $s\in \sbh_*(D)$ с мерой Рисса $\nu_s:=\nu_u-n_{\sf Z}\in \Meas^+(D)$, для которой 
\begin{equation}\label{sf}
u=\ln |f_{\sf Z}|+s\leq  M\quad \text{на $D$}.
\end{equation}  
\begin{lemma}\label{lu0s} Пусть $u_0,s,M\in \sbh_*(D)$ и имеет место неравенство 
\begin{equation}\label{u0s}
u_0+s\leq M\quad \text{на $D$.}
\end{equation} 
Тогда найдется функция  $g\in \Hol_* (D)$, с которой выполнено неравенство 
\begin{equation}\label{u0g}
u_0+\ln|g|\leq  M^{\uparrow}\quad \text{на $D$},
\end{equation} 
где поднятие $M^{\uparrow}$ определено в п.~{\rm \ref{ssM}}  в зависимости от вида области $D$  по пунктам  \eqref{Di}--\eqref{Diii}, исходя из  произвольного выбора непрерывной функции $r\colon D\to (0,1]$, удовлетворяющей  \eqref{rdD}, в \eqref{Dg}, \eqref{df:Nd}, а также чисел $\e>0$ 
в \eqref{Dg} и $P>0$ в \eqref{df:Nd_C}. 
\end{lemma}
\begin{proof}[леммы \ref{lu0s}] \textit{Случай \eqref{Di}}.\/ В силу субгармоничности функции $u_0$, усредняя  по кругам $D(z,r)$ по мере Лебега $\lambda$ обе части неравенства \eqref{u0s}, получаем
	\begin{equation}\label{usM}
u_0+s^{*r}\leq u_0^{*r}+s^{*r}\leq M^{*r}\quad \text{на $D$}.
	\end{equation}
По \cite[теорема 3]{KhaBai16} существует функция $g\in \Hol_*(D)$, для которой
\begin{equation}\label{gse}
\ln \bigl|g(z)\bigr|\leq s^{*r}(z)+\ln\frac{1}{r(z)}+(1+\e)\ln (1+|z|), \quad z\in D,
\end{equation}
откуда согласно \eqref{usM} и определению \eqref{Dg} получаем \eqref{u0g}.

\textit{Случай \eqref{Dii}}.\/  Ситуация с $\CC_{\infty}\setminus \clos D\neq \varnothing$  содержится в \cite[теорема 1]{KhaKha16}  и частично в \cite[теорема 1]{Kha15}. Для односвязной в $\CC_{\infty}$ области $D\subset \CC$  доказательство проходит по той же схеме, что и в предшествующем случае, но с  
использованием вместо \eqref{gse} неравенства 
\begin{equation*}
\ln \bigl|g(z)\bigr|\leq s^{*r}(z)+\ln\frac{1}{r(z)},\quad z\in D,
\end{equation*} 
основанного на \cite[следствие 3(iii)]{KhaBai16}. 

 \textit{Случай \eqref{Diii}}. Ситуация с $D=\CC$ разобрана в \cite[теорема 1]{KhaKha16} и частично в \cite[теорема 1]{Kha15}.
\end{proof}

По лемме \ref{lu0s} из неравенства \eqref{sf}, записываемого в виде \eqref{u0s} при $u_0:=\ln |f|$, получаем заключение \eqref{u0g}, которое означает, что $\ln |f_{\sf Z}g|=\ln|f_{\sf Z}|+\ln |g|\leq M^{\uparrow}$.  Таким образом функция $f:=f_{\sf Z}g\in \Hol_*(D)$, обращающаяся в нуль на ${\sf Z}$, искомая. 

\section{Обращение с функциями Грина}\label{Grf} 

В теореме обращения этого \S~3 используются лишь продолженные нулем функции Грина 
\cite{HK} специальной системы  относительно компактных регулярных подобластей в  $D$, содержащих некоторую фиксированную подобласть $D_0\Subset D$ с фиксированным полюсом $z_0\in D_0$. Отметим, что каждая такая {\it функция Грина --- тестовая субгармоническая финитная функция для области\/ $D$ вне подобласти\/ $D_0$}. Здесь, в отличие от теорем \ref{th1}--\ref{th:inver}, {\it область\/ $D\subset \CC_{\infty}$ произвольная.} 
Всюду в \S~\ref{Grf} в дополнение к \eqref{D} предполагаем, что  $z_0\in D_0\Subset D\subset \CC_{\infty}\neq D$.

\begin{definition}[{\rm (см. \cite[определение 1]{Kh07}, \cite[определение 1]{KhS09})}]\label{opexd} 
 Систему {\it регулярных\/} для задачи Дирихле областей $\mathcal U_{D_0} (D)\subset \{D'\Subset D \colon 	\; D_0\subset D' \}$ называем {\it  регулярной оптимально исчерпывающей в\/} $D$ с центром $D_0$, если 
$\bigcup \bigl\{D'\colon D' \in \mathcal U_{D_0} (D)\bigr\}=D$
и для любых областей $D_1$ и $D_2$, удовлетворяющих включениям  $D_0\subset D_1\Subset D_2\subset D$, выполнены два    условия:
\begin{enumerate}[{1)}]
	\item  {\it  можно подобрать область\/ $D'\in \mathcal U_{D'}(D )$ так, что\/  $D_1\Subset D' \Subset D_2$ и каждая непустая ограниченная компонента связности множества\/ ${\mathbb C}_{\infty} \setminus D'$ пересекает\/ ${\mathbb C}_{\infty} \setminus D_2$},
\item
 {\it для любой области $D'\in  \mathcal U_{D_0}(D )$ найдется такая область $D'' \in  \mathcal U_{D_0}(D )$, что имеют место включения $D_1\Subset D'' \Subset D_2$  и объединение 
$D'' \cup D'$ также принадлежит\/ $\mathcal U_{D_0}(D )$},
\end{enumerate}
и, кроме того,  эта система $\mathcal U_{D_0} (D)$ {\it условно инвариантна относительно сдвига в\/} $D$, т.\,е. из условия $D' \in \mathcal U_{D_0} (D )$, $z\in {\mathbb C}$ и $D_0\subset D'+z \Subset D$ следует, что $D'+z \in \mathcal U_{D_0} (D )$. 
\end{definition}

\begin{example} Простым  примером  регулярной оптимально исчерпывающей системы областей 
 может служить {\it специальная система всевозможных связных объединений $D'\supset D_0$ конечного числа кругов $D(z,t)\Subset D$, исключая те области $D'$, в дополнении ${\mathbb C}_{\infty} \setminus D'$ которых есть изолированные  точки}.   С такими же исключениями круги в этом примере можно заменить на относительно компактные в $D$ всевозможные $n$-угольники или, более общ\'о, односвязные подобласти \cite[{теоремы} 4.2.1, 4.2.2]{Rans}, \cite[{2.6.3}]{HK} какого-либо специального вида.  
\end{example}

\begin{theorem}\label{th3} Пусть   $M=M_+-M_-\in \dsbh_* (D)$  с зарядом Рисса $\nu_M\in \Meas(D)$, где $M_+\in \sbh_* (D)\cap C(D)$ и $M_-\in \sbh_*(D)$, а также $z_0\in D_0\cap \dom M \Subset D$.
Пусть для меры $\nu\in \Meas^+(D)$  при некотором  $r_0\in \RR_*^+$ выполнено   \eqref{dMa+}. Пусть\/ $\mathcal U_{D_0} (D)$ --- 	  {\it  регулярная оптимально исчерпывающая  система областей в\/} $D$ с центром $D_0$, для которой    с некоторой постоянной $C\in \RR$  выполнены неравенства
\begin{equation}\label{estv+g}
\int\limits_{ D \setminus \{z_0\}}  g_{D'}(\cdot ,z_0) \dd {\nu} 		\leq	\int\limits_{D\setminus \{z_0\}} g_{D'}(\cdot ,z_0) \dd {\nu}_M	+C \quad \text{для всех $D'\in \mathcal U_{D_0} (D)$},
 \end{equation}
т.\,е. заряд $\nu_M$  --- аффинное выметание меры  $\nu$ для $D$ вне ${z_0}$ относительно класса  функций Грина $g_{D'}(\cdot ,z_0)$ с  $D'\in \mathcal U_{D_0} (D)$. 
Тогда  найдется   функция $u\in \sbh_*(D)$ с мерой Рисса $\nu_u\geq \nu$,  удовлетворяющая неравенству   $u\leq M$ на $D$.
\end{theorem}
\begin{proof}  Пусть  $\nu_{M_+}$ и $\nu_{M_-}$ --- меры Рисса соответственно  функций ${M_+}$ и ${M_-}$. Тогда равномерную по постоянной  $C$ серию неравенств  \eqref{estv+g}   можно записать в обозначении $\nu_1:=\nu +\nu_{M_-}$ в виде
\begin{equation}\label{estv+g+}
\int\limits_{ D \setminus \{z_0\}}  g_{D'}(\cdot ,z_0) \dd {\nu}_1 		\leq	\int\limits_{D\setminus \{z_0\}} g_{D'}(\cdot ,z_0) \dd \nu_{M_+}	+C
\text{ для всех $D'\in \mathcal U_{D_0} (D)$},
 \end{equation}
где $\nu_1, \nu_{M_+}\in \Meas^+(D)$ --- уже {\it положительные меры.\/}  
Выберем какую-нибудь субгармоническую в $D$ функцию $u_1\in \sbh_*(D)$ с мерой Рисса $\nu_1$.  Для ее меры Рисса $\nu_1$,  
ввиду  условия $z_0\in \dom M $, а также  условия \eqref{dMa+} на $z_0$, выполнено условие  \eqref{dMa+} с  заменой $\nu$ на $\nu_1$. Следовательно,  $M_-(z_0)\neq -\infty$ и  обязательно $u_1(z_0)\neq -\infty$.  Далее потребуются 
 вариации утверждений из \cite[основная теорема, теорема 6]{Kh07}:

\begin{theoB}[{\rm (частный случай  \cite[теорема (основная)]{KhS09})}]
Пусть функция\/  $M \in \sbh_*(D)$  с мерой Рисса $\nu_{M}$ ограничена снизу в некоторой открытой окрестности замыкания   $\clos D_0$,   $u\in \sbh_*(D)$ --- функция с мерой Рисса $\nu$ на $D$ и $u(z_0)\neq -\infty$,  система областей\/ $\mathcal U_{D_0} (D)$ --- 	регулярная оптимально исчерпывающая для $D$ с центром $D_0\ni z_0$. 
Если\footnote{К сожалению,  в формулировке основной теоремы из нашей работы \cite{Kh07}, на промежуточном этапе доказательства которой и основано  \cite[теорема (основная)]{KhS09},  допущена досадная опечатка в знаках $\pm$. Так, используемое в ее формулировке соотношение  
\cite[п.~(h1), (2.11)]{Kh07} должно выглядеть в точности как \eqref{est:g}.  Дальнейший комментарий --- в сноске к \cite[теорема (основная)]{KhS09}.} 
\begin{equation}\label{est:g}
-\infty < \inf_{D' \in \mathcal U_{D_0}(D )} \left(-\int_{D \setminus \{z_0\}} g_{D'}(\cdot  , z_0) \dd\nu_u+
\int_{D\setminus \{z_0\}} g_{D'} (\cdot , z_0 )\dd\nu_M\right),   
\end{equation}  
то  для любой непрерывной функции\/  $r \colon D\to \RR^+$, удовлетворяющей условию\/ 
\eqref{rdD},  найдется функция\/ $v\in \sbh_*(D)$, гармоническая в открытой окрестности точки $z_0$, для которой   $u +v\overset{\eqref{M*r}}{\leqslant} M^{*r}$ на $D$. При этом  если еще и  $M\in C(D)$, то переменное усреднение $M^{*r}$ в правой части последнего неравенства  можно заменить на $M$.
\end{theoB}
 По теореме B, примененной  к функциям $u_1$ и {\it непрерывной\/} функции 
$M_+$ вместо соответственно $u$ и $M$, ввиду \eqref{estv+g+}, соответствующего условию \eqref{est:g}, найдется функция $v \in \sbh_* (D)$, гармоническая в окрестности точки $z_0$, с которой $u_1+v\leq M_+$ на $D$. По построению   $u_1\in \sbh_* (D)$ с  мерой Рисса  $\nu_1:=\nu +\nu_{M_-}$. Следовательно,  мера Рисса функции $u_0:=u_1-M_-$ --- это мера $\nu$, т.\,е. существует функция $u:=u_0+v \in \sbh_*(D)$ с мерой Рисса $\nu_u\geq \nu$, для которой $u\leq M_+-M_-=M$ на $D$, что и завершает доказательство теоремы \ref{th3}.
\end{proof}

\begin{corollary}\label{holDg} Пусть в условиях теоремы\/ {\rm \ref{th3}}   область $D\subset \CC$, 
$M\in \sbh_*(D)\cap C(D)$ и  мера $\nu_M$  --- аффинное выметание последовательности  ${\mathsf Z}=\{{\mathsf z}_k\}_{k=1,2,\dots}\subset D$, $z_0\notin {\sf Z}$,  для $D$ вне $S:=\{z_0\}$ относительно класса  функций Грина $g_{D'}(\cdot ,z_0)$ с  $D'\in \mathcal U_{D_0} (D)$,  
т.\,е. для некоторого числа $C\in \RR$ выполнено условие \eqref{estv+g} в виде 
\begin{equation*}
\sum\limits_{ {\mathsf z}_k\in D' }  g_{D'}({\mathsf z}_k ,z_0)  \overset{\eqref{nufZ}, \eqref{balZ}}{\leq}	\int\limits_{D\setminus \{z_0\}} g_{D'}(\cdot ,z_0) \dd {\nu}_M	+C \quad \text{ для всех $D'\in \mathcal U_{D_0} (D)$}.
 \end{equation*}
Тогда   найдется  функция $f\in \Hol_* (D)$, для которой  $ f({\mathsf Z})=0$ и выполнено  неравенство $|f|\leq \exp M^{\uparrow}$ на $D$, где поднятие $M^{\uparrow}$ определено в п.~{\rm \ref{ssM}} в зависимости от вида области $D$  по пунктам  \eqref{Di}--\eqref{Diii}, исходя из  произвольного выбора непрерывной функции $r\colon D\to (0,1]$, удовлетворяющей  \eqref{rdD}, в \eqref{Dg}, \eqref{df:Nd}, а также чисел $\e>0$ в \eqref{Dg} и $P>0$ в \eqref{df:Nd_C}.  
\end{corollary}
Выводится из теоремы  \ref{th3} так же, как и импликация {\bf h\ref{h3}}$\,\Rightarrow\,${\bf h\ref{h1}} с поднятием $M^{\uparrow}$ из  теоремы \ref{th:inver} в подразделе \ref{h31}.  

\begin{remark} Регулярную  оптимально исчерпывающую   систему областей $\mathcal U_{D_0} (D)$ с центром $D_0\subset D$ в теореме \ref{th3} и следствии \ref{holDg} на основе анализа тонких совместных результатов В.~Хансена и И.~Нетуки \cite{HN} об аппроксимации мер Йенсена гармоническими мерами можно заменить на систему областей $D'\Subset D$, включающих в себя  область $D_0\Subset D$ и полученных из  исчерпывающей область $D$ последовательности регулярных областей $D_n\Subset D$, $n\in \NN$  с аналитической или кусочно линейной или иной <<хорошей>> границей удалением из областей $D_n$ всевозможных произвольных  конечных  наборов  попарно не пересекающихся замкнутых кругов. При этом для полученной системы областей надо все-таки требовать условную инвариантность относительно сдвига в $D$ из определения \ref{opexd}.
\end{remark}

\section{Обращение с аналитическими и полиномиальными дисками}\label{apnd}

Важный подкласс в классе $J_{z_0}(D)$ мер Йенсена порождают аналитические диски в $D$ с центром $z_0$.
{\it Аналитическим замкнутым диском  в области\/ $D$ с центром\/}  $z_0\in D$\/ называется функция $g \colon  \clos \DD \to D$, непрерывная на $\clos \DD$ с  голоморфным сужением в $\DD$, для которой $g(0)=z_0$ \cite{BSch}--\cite{CR+}. В частности, $g(\clos \DD)\Subset D$. Для любого такого аналитического замкнутого диска $g$ легко показать, что функция от $w\in \CC_{\infty}\setminus \{z_0\}$, равная 
\begin{subequations}\label{df:Vmu+}
\begin{align}
&\frac{1}{2\pi} \int_0^{2\pi} \ln \Bigl|\frac{w-g(e^{i\theta})}{w-z_0}\Bigr| \dd \theta
= \frac{1}{2\pi} \int_0^{2\pi} \ln \Bigl|1-\frac{g(e^{i\theta})-z_0}{w-z_0}\Bigr| \dd \theta \quad\text{\it при $z_0\neq \infty$},
\tag{\ref{df:Vmu+}o}\label{{df:Vmu}o+}
\\
\intertext{где при $w=\infty$ подынтегральные выражения доопределены значением  $0$ (ср. с \eqref{{df:Vmu}o}),} 
&\frac{1}{2\pi} \int_0^{2\pi} \ln \Bigl|\frac{w-g(e^{i\theta})}{g(e^{i\theta})}\Bigr| \dd \theta=
\frac{1}{2\pi} \int_0^{2\pi} \ln \Bigl|1-\frac{w}{g(e^{i\theta})}\Bigr| \dd \theta \quad\text{\it при $z_0= \infty$},
\tag{\ref{df:Vmu+}$\infty$}\label{{df:Vmu}in+}
\end{align}
\end{subequations}
где при $g(e^{i\theta})=\infty$ подынтегральные выражения доопределены значением  $0$ (ср. с \eqref{{df:Vmu}in}),  --- потенциал Йенсена внутри области $D$ с полюсом в точке $z_0$.  В частности,   \eqref{df:Vmu+} определяет тестовую субгармоническую положительную финитную функцию для $D$ вне $\{z_0\}$. В замечании \ref{r2} мы называли функции \eqref{df:Vmu+}
 логарифмическими потенциалами аналитических дисков. Если аналитический диск  $g$ в области $D$ с центром $z_0\in D$ --- {\it многочлен\/} комплексной переменной,  то естественно называть его {\it полиномиальным диском в $D$ с центром в\/}  $z_0\in D$. 
\begin{theorem}\label{th:4} Пусть функция $M\in \dsbh_* (D)$ с зарядом Рисса $\nu_M$,  точка $z_0\neq \infty$ и мера $\nu\in \Meas^+(D)$ такие же, как в теореме\/ {\rm \ref{th3}}. Если заряд $\nu_M$ ---  аффинное выметание меры $\nu$ для $D$ вне $S:=\{z_0\}$ относительно класса функций \eqref{{df:Vmu}o+}, т.\,е. существует постоянная  $C\in \RR$, с которой выполнено  неравенство 
\begin{equation*}
\int_{D} \frac{1}{2\pi}\int_0^{2\pi} \ln \left|1-\frac{g(e^{i\theta})-z_0}{z-z_0}\right| \dd\theta \dd \nu (z)	
\overset{\eqref{{df:Vmu}o+}}{\leq}  	\int_{D}  \frac{1}{2\pi}\int_0^{2\pi} \ln \left|1-\frac{g(e^{i\theta})-z_0}{z-z_0}\right| \dd\theta \dd \nu_M(z)+C 	
\end{equation*}
для всех аналитических замкнутых  или лишь полиномиальных  дисков $g$ в $D$ с центром $z_0$, то  
найдется такая  функция $u\in \sbh_*(D)$ с мерой Рисса $\nu_u\geq \nu$, что  $u\leq M$ на $D$.
\end{theorem}
В случае субгармонической  функции $M$ обсуждение схемы доказательства теоремы  \ref{th:4} содержится  в  \cite[1.2.1--1.2.2, дополнения 1.2.3, 1.2.4]{Khsur}. Это одна из причин, по которой  мы опускаем здесь 
доказательство  теоремы  \ref{th:4}. Другая в том, что многомерный вариант теоремы \ref{th:4} в $\CC^n$ более естественен и будет рассмотрен с применениями в ином месте. 

Как импликация {\bf h\ref{h3}}$\,\Rightarrow\,${\bf h\ref{h1}} с поднятием $M^{\uparrow}$ в подразделе \ref{h31} и следствие \ref{holDg} из теоремы \ref{th:4} выводится 
\begin{corollary}\label{holDgad} В условиях теоремы\/ {\rm \ref{th:4}}  вместо меры $\nu$
рассмотрим последовательность точек\/  ${\mathsf Z}=\{{\mathsf z}_k\}_{k=1,2,\dots}\subset D\subset \CC$, $z_0\in D\setminus {\sf Z}$ и предполагаем, что $M\in \sbh_*(D)\cap C(D)$. Если мера $\nu_M$ ---  аффинное выметание  последовательности\/  ${\mathsf Z}$ для $D$ вне $S:=\{z_0\}$ относительно класса функций \eqref{{df:Vmu}o+}, т.\,е. существует постоянная  $C\in \RR$, с которой неравенство 
\begin{equation*}
	 \sum_{ {\mathsf z}_k\in D} \int_0^{2\pi}  \log \left|1-\frac{g(e^{i\theta})-z_0}{{\mathsf z}_k-z_0}\right| \dd\theta  	
	 \leq  	\int_{D}  \int_0^{2\pi} \log \left|1-\frac{g(e^{i\theta})-z_0}{z-z_0}\right| \dd\theta \dd\nu_M(z)+C 	
\end{equation*}
выполнено для всех аналитических замкнутых или только полиномиальных  дисков $g$ в $D$ с центром $z_0$, то  найдется функция $f\in \Hol_* (D)$, для которой  $ f({\mathsf Z})=0$ и  $|f|\overset{\ref{ssM}}{\leq} \exp  M^{\uparrow}$ на $D$, где поднятие $M^{\uparrow}$ определено в п.~{\rm \ref{ssM}} с комментарием из заключения следствия\/ {\rm \ref{holDg}}.
\end{corollary}


\begin{thebibliography}{86}

\RBibitem{KhaRoz18}
\by Б.~Н.~Хабибуллин, А.~П.~Розит
\paper К распределению нулевых множеств голоморфных функций
\jour Функц. анализ и его прил.
\yr 2018
\vol 52
\issue 1
\pages 26--42

\RBibitem{KhaKha19}
\by Б.~Н.~Хабибуллин, Э.~Б.~Хабибуллина
\paper К распределению нулевых множеств голоморфных функций. II
\paperinfo краткое сообщение
\jour Функц. анализ и его прил.
\yr 2018
%%\vol 
%%\issue 
%%\pages 
\toappear

\RBibitem{KhT16} 
\by Б.\,Н.~Хабибуллин, Н.\,Р.~Таминдарова
\paper Распределение нулей и масс голоморфных и субгармонических функций:
условия типа Адамара и Бляшке
\yr 2018
\finalinfo  \href{https://arxiv.org/abs/1512.04610v4}{https://arxiv.org/abs/1512.04610v4}\;

\Bibitem{Rans}
\by Th.~Ransford
\book Potential Theory in the Complex Plane
\publ Cambridge University Press
\publaddr Cambridge
\yr 1995

\RBibitem{L} 
\by Н.\,С.~Ландкоф 
\book Основы современной теории потенциала 
\publ Наука 
\publaddr М.
\yr 1966

\RBibitem{Kh07}
\by Б.~Н.~Хабибуллин
\paper Последовательности нулей голоморфных функций,
представление мероморфных функций и гармонические миноранты
\jour Матем. сб.
\yr 2007
\vol 198
\issue 2
\pages 121--160

\RBibitem{Khsur} 
\by Б.~Н.~Хабибуллин
\book  Полнота систем экспонент и множества единственности  
\publ    РИЦ БашГУ 
\publaddr  Уфа 
\bookinfo издание четвёртое 
%% (1-ое --- 2006; 2-ое --- 2008; 3-е --- 2011), дополненное%%, 
%%\isbn  %%ISBN 978-5-7477-2540-9 
 %%\totalpages 192 
\yr 2012
\finalinfo  \href{http://www.mathnet.ru/rus/person/8650}{http://www.mathnet.ru/rus/person/8650}\;
%%http://www.researchgate.net/profile/Bulat\_Khabibullin
%%https://www.researchgate.net/publication/271841461_Polnota_sistem_eksponent_i_mnozestva_edinst%%ven%%nosti_Completeness_of_Exponential_Systems_and_Uniqueness_Sets


\Bibitem{BH}
\by J.~Bliedtner, W.~Hansen 
\book Potential Theory --- An Analytic and Probabilistic Approach to Balayage
\yr 1986
\publ Springer-Verlag
\publaddr Berlin--Heidelberg--N.Y.--Tokyo

\RBibitem{B} 
\by Н.~Бурбаки 
\book Интегрирование. Меры, интегрирование мер 
\publ Наука 
\publaddr М. 
\yr 1967

\RBibitem{HK} 
\by У.~Хейман, П.~Кеннеди
\book Субгармонические функции
\publ Мир
\publaddr М.
\yr 1980

\RBibitem{Kha15}
\by Б.\,Н.~Хабибуллин
\paper  Последовательности неединственности для весовых пространств голоморфных функций 
\jour Изв. вузов. Матем.
\yr 2015
\vol {\rm 4}
%%\issue 4
\pages 75--84


\RBibitem{KhaKha16}
\by Б.~Н.~Хабибуллин, Ф.~Б.~Хабибуллин
\paper О множествах неединственности для пространств голоморфных функций
\jour Вестн. Волгогр. гос. ун-та. Сер. 1, Мат. Физ.
\yr 2016
\vol {\rm 4(35)}
\pages 108--115

\RBibitem{Khab03}
\by Б.~Н.~Хабибуллин
\paper Критерии (суб-)гармоничности и~продолжение (суб-)гармонических функций
\jour Сиб. матем. журн.
\yr 2003
\vol 44
\issue 4
\pages 905--925

\RBibitem{KhaBai16}
\by Б.~Н.~Хабибуллин, Т.~Ю.~Байгускаров
\paper Логарифм модуля голоморфной функции как миноранта для субгармонической функции
\jour Матем. заметки
\yr 2016
\vol 99
\issue 4
\pages 588--602

\RBibitem{KhS09}
\by  Б.\,Н.~Хабибуллин
\paper Нули голоморфных функций с ограничениями на рост в области
\inbook  Исследования по математическому анализу
\serial Математический форум (Итоги науки. Юг России)
\vol 3
\publ  Владикавказский научный центр РАН и РСО--А
\publaddr  Владикавказ  
\pages  282--291
\yr 2009 
\finalinfo \href{http://elibrary.ru/download/90353677.pdf}{http://elibrary.ru/download/90353677.pdf} \;

\bibitem{HN} 
\by W.~Hansen, I.~Netuka
\paper Convexity properties of harmonic measures
\jour Adv. Math. 
\yr 2008
\vol 218
\issue 4
\pages 1181--1223


\Bibitem{BSch}
S.~Bu, W.~Schachermayer W.
\paper  Approximation of Jensen me\-a\-su\-res by image   me\-a\-su\-res  under holomorphic functions and applications 
\jour Trans. Amer. Math. Soc. 
\yr 1992
\vol 331
\issue 2
\pages 585--608

\Bibitem{Poletsky49} 
\by E.\,A.~Poletsky
\paper Holomorphic currents
\jour Indiana Univ. Math. J.
\yr 1993
\vol 42
%%\issue 4
\pages 85--144

\Bibitem{Po99}
\by E.\,A.~Poletsky
\paper Disk envelopes of functions II
\jour J. Funct. Anal.
\yr 1999
\vol 163
%%\issue 4
\pages 111--132

\Bibitem{CR+}
\by B.\,J.~Cole, T.\,J.~Ransford
\paper Jensen measures and harmonic measures
\jour  J. Reine Angew. Math. 
\yr  2001
\vol 541
%%\issue 2
\pages  29--53
\end{thebibliography}
\end{document}